\documentclass[10pt]{article}
\usepackage{amssymb}
\usepackage{enumerate}
\usepackage{graphicx}

\def\R{{\rm I\! R}}

\newtheorem{theorem}{Theorem}
\newtheorem{corollary}{Corollary}

\title{Asymptotic stability and periodic solutions  of vector Li\'enard equations}
  
\author{F. Briata \footnote{Dip. di Scienze, 
Universit\`a degli Studi \lq\lq G.d'Annunzio\rq\rq
, V.le Pindaro, 42, 65127 Pescara - Italy. 
Email: briata@sci.unich.it,
Phone: ++39(085)4537692} , 
M. Sabatini  \footnote{Corresponding author - Dip. di Matematica, Univ. di Trento, I-38050 Povo, (TN) - Italy.
Email: marco.sabatini@unitn.it,
Phone: ++39(0461)881670, Fax: ++39(0461)881624}
}
\begin{document}
\maketitle
\begin{abstract} We prove the asymptotic stability of the equilibrium solution of a class of vector Li\'enard equations by means of LaSalle invariance principle. The key hypothesis consists in assuming that the intersections of the manifolds in $\{\dot V = 0\}$ be isolated. We deduce an existence theorem for periodic solutions of periodically perturbed vector Li\'enard equations.


\end{abstract}

\section{Introduction}

In this paper we are concerned with the vector Li\'enard equation
\begin{equation} \label{equaliev}
X'' + f(X) \cdot X' + g(X) = 0,
\end{equation}
where $f \in C^0(\Omega,\R^n)$, $\Omega $ open connected subset of $\R^n$, $g\in Lip(\prod_{i=1}^n (a_i,b_i),\R^n)$, the space of locally lipschitzian functions on $\prod_{i=1}^n (a_i,b_i)$, and $\cdot$ denotes the inner product of $\R^n$. We assume the origin $O$ to be an isolated equilibrium point of the equivalent system
\begin{equation} \label{sysliev}
X' = Y, \qquad Y' =  - g(X) - f(X) \cdot  Y,
\end{equation}
In the planar case ($n=1$), the equation (\ref{equaliev}) has been widely studied from many points of view: stability of equilibria, boundedness and oscillation of solutions, existence, uniqueness or multiplicity of periodic solutions (see \cite{P},\cite{SC},\cite{V}). The higher dimensional case present substantial additional difficulties, which prevent straightforward extensions of planar results. This is in particular true for an equilibrium's stability properties, which are essential in studying the dynamics of perturbed sytems. Even if such properties are often analyzed by linearizing the system at the equilibrium, this cannot be done if some component of $f(X)$ or some partial derivative of $g(X)$ vanish at such a critical point.

Among the few available results, we point out two papers. The first one is concerned with the study  of the singular perturbation problem $\varepsilon\ddot x+R(x)\dot x+g(x)=0$ for $0<\varepsilon\ll 1$,  performed in \cite{PW}. For small positive values of $\varepsilon $ the solutions approach the subspace $\dot x=-R^{-1}(x)g(x)$. 

A class of Li\'enard vector equations was considered in \cite{A}, under quite restrictive hypotheses. In \cite{A} $g(X)$ is assumed to be the gradient of a suitable positive definite homogeneous scalar function, $f(X)$ to be the jacobian matrix of a homogeneous vector field $F(X)$ on $\R^n$, such that $F(X) \cdot g(X)$ be positive definite. This ensures that the conditions required by Barbashin-Krasokskii-LaSalle theorem (\cite{V}, thm. 79) hold for the system
\begin{equation} \label{sysalex}
X' = Y - F(X), \qquad Y' =  - g(X) ,
\end{equation}
which is equivalent to (\ref{equaliev}). In fact, such hypotheses  are satisifed by taking  the energy  $V(X,Y) =  |Y|^2/2+ G(X)$ as a Liapunov function, since the derivative $\dot V$ of $V$ along the  solutions of (\ref{sysalex}) vanishes only on the set $\{X = 0\}$, i. e. on the \lq\lq Y axis\rq\rq. The non-positive invariance of such a set is an immediate consequence of $X' = Y$ on the on the set $\{X = 0\}$. 

In this paper we study the asymptotic stability of the origin for the system (\ref{sysliev}) under different hypotheses, admitting the possibility of  a larger vanishing set $\{\dot V = 0\}$, possibily containing several hypersurfaces not contained in the set $\{X = 0\}$. This allows to apply our result to equations excluded by the approach of \cite{A}, as
\begin{equation}  \label{esempio1} 
\left\{\begin{array}{cl}
x_1 '  &=  y _1   \\ 
x_2 '  &=  y _2   \\ 
\end{array}\right.   \qquad \qquad
\left\{\begin{array}{cl}
y_1 '  &= - x_1  - y_1 (x_1 - x_2)^2 \\
y_2 '  &= - x_2  - y_2 (x_1 + x_2^2)^2. 
\end{array}\right.  
\end{equation}
We assume $g(X) = (g_1(x_1), \dots, g_n(x_n))$, so that the coupling of the $n$ scalar equations in (\ref{equaliev}) is entirely due to the dissipative terms $ f(X) \cdot X'$.  
Then the system (\ref{sysliev}) appears as a Li\'enard-like perturbation of a natural Hamiltonian system with separable variables. Actually, the above Liapunov function $V$ is the hamiltonian function of the unperturbed system. We prove the  equilibrium's asymptotic stability assuming that the intersection of the manifolds in $\dot V = 0$ consists of isolated points (theorem \ref{asystab}). In order to do this, we prove that every solutions starting at a point of any of such manifolds leaves it immediately, proving the non-invariance of $\dot V = 0$. In corollary \ref{stabglob} we extend our result to prove the equilibrium's global asymptotic stability. In corollary \ref{persol}, we apply theorem  \ref{asystab} in order to prove the existence of periodic solutions for small perturbations of the type
\begin{equation} \label{equalievper}
X'' + f(X) \cdot X' + g(X) = h(t,X,X',\varepsilon),
\end{equation}
where $h(t,X,X',\varepsilon) \equiv h(t+T,X,X',\varepsilon) \in C^1(\R^{2n+1} \times [0,\overline{\varepsilon}),\R^n)$, $T>0$.

\section{Li\'enard equation, $n > 1$} \label{Lienard ordn}

Let us set $X = (x_1, \dots, x_n)$, $Y = (y_1, \dots,y_n)$,  $Z = (x_1,\dots,x_n,y_1 \dots, ,y_n)$. Let us consider a map $f=(f_1,\dots,f_n)$ defined on an open connected subset $\Omega$ of $\R^n$ containing the $n$-origin $O_X$. For a scalar map $g:(a,b) \rightarrow \R$, we say that $g \in Lip((a,b),\R)$ if $g$ is locally lipschitzian in all of $(a,b)$.

Let us consider a system of $n$ coupled Li\'enard equations %
\begin{equation} \label{equaLienardn}
\left\{\begin{array}{cl}
x_1 '' +f_1(x_1, \ldots, x_n)x_1'+g_1(x_1)  &= 0 \\ 
x_2 '' +f_2(x_1, \ldots, x_n)x_2'+g_2(x_2)  &= 0 \\
\vdots \\
x_n '' +f_n(x_1, \ldots, x_n)x_n'+g_n(x_n)  &= 0, 
\end{array}\right.  
\end{equation}
where $f \in C^0(\Omega,\R)$, $g_i(x_i) \in Lip((a_i,b_i), \R)$, $a_i < 0 < b_i$, $i=1,\dots,n$.
The equivalent first-order system is\begin{equation}  \label{sysn} 
\left\{\begin{array}{cl}
x_1 '  &=  y _1   \\ 
\vdots\\
x_n '  &=  y _n   \\ 
\end{array}\right.   \qquad \qquad
\left\{\begin{array}{cl}
y_1 '  &= -g_1(x_1) -y_1 f_1(x_1, \ldots, x_n) \\
\vdots\\
y_n '  &= -g_n(x_n) -y_n f_n(x_1, \ldots, x_n) . 
\end{array}\right.  
\end{equation}
We assume $g_i(0) =0$, $x_ig_i(x_i) > 0$ for $x_i \in (a_i ,0) \cup(0, b_i)$, $i= 1,\dots,n$, so that the $2n$-origin $O = (O_X,O_Y)$ is an isolated equilibrium point of (\ref{sysn}). Also in this case, if the $g_i$'s are differentiable at $0$, in order to study the stability of $O$ one could compute the eigenvalues of the linearized system at $O$, which is non-degenerate if and only if $\prod_{i=1}^n g_i'(0)  \neq 0$. Moreover, the eigenvalues' real parts depend on the values of $f_i(O_X)$, $i=1,\dots,n$, similarly to what  happens for $n=1$,
$$
\lambda^{1,2}_i = \frac{-f_i(O_X) \pm \sqrt{f_i(O_X)^2 - 4g_i'(0)}}{2},
$$
Hence, if   $f_i(O_X)$ vanishes and $g_i'(0) >0$ for some $1 \leq i \leq n$, a different approach is required. This occurs also if $g_i'(0)=0$, for some $1 \leq i \leq n$. 
 
The Liapunov function approach is more general, in two ways. It applies to degenerate critical points, and allows also to give an estimate of the equilibrium's region of attraction. 
In order to apply LaSalle invariance principle to a suitable Liapunov function $V$, we have to prove that the vanishing set $\{\dot V = 0\}$ is not positively invariant. Such a set is the set-theoretical union of  hypersurfaces. In next theorem we prove that, if  such hypersurfaces' intersections are isolated, then $\{\dot V = 0\}$ is not positively invariant.  
 From now on we use  the word {\it isolated} referring to  the topology of $\R^n$. For definitions abopu dynamical systems, we refer to \cite{BS}.
 
\begin{theorem} \label{asystab}  If  
\begin{itemize}
\item[1)] $x_ig_i(x_i) > 0$ for $x_i \in  (a_i,b_i)$,   $x_i \neq 0$, $i= 1,\dots, n$;
\item[2)] $f_i(X) \geq 0$, in $\Omega$, $i= 1,\dots, n$; 
\item[3)] for every subset $S \subset \{ 1,\dots,n \}$  the set $\{ x_i =0, {\rm\ for\ } i \in S\} \cap \{ f_i =0, {\rm\ for\ } i \not \in S\}$ consists of isolated points;
\item[4)]  $\cap_{i=1}^n \{ f_i = 0 \}$ consists of isolated points;
\end{itemize}
then $O$ is asymptotically stable.
\end{theorem}
{\it Proof.}
Let us set
\begin{equation} 
V(x_1,\dots, x_n,y_1,\dots,y_n) = \sum_{i=1}^n \left( G_i(x_i)+ \frac{y_i^2}{2} \right),
\end{equation}
 where $G_i(x)=  \int_{0}^{x_i}g_i(s)ds$, with $i=1,\dots, n$. By the above assumptions, $V$ is positive definite at $O$.
In order to study the stability of the equilibrium point $O$, we study the sign of the derivative of $V$ along the solutions of ( \ref{sysn} )
$$
\dot{V}(Z) = -  \sum_{i=1}^n y_i^2 f_i(X)   .
$$
The hypothesis 2) implies that $\dot{V}(Z) \leq 0$ in a neighbourhood of $O$, hence the origin is stable. 
In order to apply the LaSalle principle, we study the orbits' behaviour on the set $ W_O =\{  \sum_{i=1}^n y_i^2 f_i(X) = 0 \}$. The equality $ \sum_{i=1}^n y_i^2 f_i(X) = 0$ holds if and only if one of the following holds: \\

$a)$ $y_1 = \dots = y_n=0$, 

$b)$ for some subset $S \subset \{ 1,\dots,n \}$ one has  $ y_i  = 0$ for $i \in S$, and $ f_i  = 0$ for $i \not \in S$, 

$c)$ $f_1 = \dots = f_n=0$. \\

Let us denote by $ W_O^a$, $ W_O^b$, $ W_O^c$, the set defined by the respective condition, so that $ W_O = W_O^a \cup W_O^b \cup  W_O^c$. We show that every orbit $\gamma(t)$ such that $\gamma(0) \in W_O$ leaves immediately $ W_O^a$, $ W_O^b$, $ W_O^c$, hence it leaves immediately $W_O$. In other words, we prove that if $\gamma(0) \in W_O$, then for all $t\in (0,\epsilon)$ one has $\gamma(t) \not\in W_O$, for some positive $\epsilon$. We also set $W_O^S = \{ y_i  = 0 {\rm \ for\ } i \in S,  f_i  = 0 {\rm\ for\ } i \not \in S\}$.

Case $a)$. 
Let $\gamma(t)$ be an orbit of (\ref{sysn}), such that $\gamma(0) \in W_O^a$.
Since $y_1(0) = 0$, for $t=0$ one has  $y_1 '(0)= -g_1(x_1(0))$ in  (\ref{sysn}).
If  $x_1 (0)\neq 0$, then  $y_1'(0) = -g_1(x_1(0)) \neq 0$, hence $\gamma$ leaves immediately the set $W_O^a$. 
As a consequence, $x_1(0) = 0$, so that $y_1 ' (0)= 0$, and at $\gamma(0)$ the system appears as follows:
\begin{equation}  
\left\{\begin{array}{cl}
x_1 '  &=  0   \\ 
x_2 '  &=  0   \\ 
\dots  \\
x_n '  &=  0   \\ 
\end{array}\right.   \qquad \qquad
\left\{\begin{array}{cl}
y_1 '  &= 0 \\
y_2 '  &= -g_2(x_2) \\
\dots  \\
y_n '  &= -g_n(x_n) . 
\end{array}\right.  
\end{equation}
We may apply to $y_i$, $i > 1$, the argument just used on $y_1$, proving that either $\gamma$ leaves immediately the set $W_O^a$, or $x_i (0)= 0$. As a consequence, one has   $y_1(0)= \dots = y_n(0)=x_1(0)= \dots =x_n(0)=0$. Hence the unique orbit that does not leave the set $y_1=y_2=0$ is the equilibrium point $O$.

Case $b)$.  Assume $\gamma(0) \in W_O^S$.
As in case a), we may prove that, if $y_i(0)=0$ for $i\in S$, then $x_i(0)=0$,  for $i\in S$, otherwise $\gamma$ leaves the set  $\{ y_i =0, {\rm\ for\ } i \in S\}$, hence it leaves immediately the set $W_O^S$. 
If $\gamma$ remains in $ W_O^S$, then it remains in the set $\{ x_i =0, {\rm\ for\ } i \in S\} \cap \{ f_i =0, {\rm\ for\ } i \not \in S\}$.
By the hypothesis 3), such a set consists of isolated points, hence there exist $\overline{X}\in \R^n$, such that $X(t) = \overline{X}$, for $t > 0$. 
As a consequence,  one has $y_i(t) = x_i'(t) = 0$,  $i=1, \dots,n$, hence $W_O^S \subset W_O^a$. Now we may use the fact that $\gamma(t)$ leaves immediately $W_O^a$, hence also  $W_O^S$. 

Case $c)$.
Let $\gamma(0)  \in  W_O^c$. By hypothesis 4), $\cap_{i=1}^n \{ f_i = 0 \}$ consists of isolated points.  If $\gamma(t)$ starts at a point of $W_O^c$ and remains in $W_O^c$, then its $x$-components are constant. Then one has $y_i (t)= x_i' (t)= 0$, for $i=1,\dots,n$, hence on $W_O^c$ the system has the following form
\begin{equation}  
\left\{\begin{array}{cl}
x_1 '  &=  0   \\ 
\dots  \\
x_n '  &=  0   \\ 
\end{array}\right.   \qquad \qquad
\left\{\begin{array}{cl}
y_1 '  &= -g_1(x_1) \\
\dots  \\
y_n '  &= -g_n(x_n) . 
\end{array}\right.  
\end{equation}
Now we can repeat the argument of case $a)$, proving that $y_1' (0)= \dots = y_n'(0) =0$, hence, at $\gamma(0)$, $x_i(0)=y_i(0) = 0$, $i=1,\dots,n$. This shows that $\gamma(0) $ remains in $W_O^c$ only if $\gamma(0) = O$.
\hfill $\clubsuit$
\bigskip

We give a couple of examples of systems satisfying the hypotheses of theorem \ref{asystab}. 
The system
\begin{equation}  \label{esempio2} 
\left\{\begin{array}{cl}
x_1 '  &=  y _1   \\ 
x_2 '  &=  y _2   \\ 
\end{array}\right.   \qquad \qquad
\left\{\begin{array}{cl}
y_1 '  &= - x_1  - y_1 x_1^2 (x_2 - 1)^2 (x_1 + 1)^2 \\
y_2 '  &= - x_2  - y_2 x_2^2 (x_1 - 1)^2 (x_2 + 1)^2. 
\end{array}\right.  
\end{equation}
One has $\{ f_1 = 0 \} \cap \{ x_2 = 0 \} =  \{ (0,0), (-1,0) \}$ and $\{ f_2 = 0 \} \cap \{ x_1 = 0 \} =  \{ (0,0), (0,-1) \}$. \\
Moreover $\{ f_1 = 0 \} \cap \{ f_2 = 0 \} =  \{ (0,0), (0,-1), (1,1), (-1,0), (-1,-1) \} $. 
\smallskip

The system
\begin{equation}  \label{esempio3} 
\left\{\begin{array}{cl}
x_1 '  &=  y _1   \\ 
x_2 '  &=  y _2   \\ 
\end{array}\right.   \qquad \qquad
\left\{\begin{array}{cl}
y_1 '  &= - x_1  - y_1 (x_1^2 + 2x_2^2 - 1)^2 \\
y_2 '  &= - x_2  - y_2  (2x_1^2 + x_2^2 - 1)^2. 
\end{array}\right.  
\end{equation}
One has $\{ f_1 = 0 \} \cap \{ x_2 = 0 \} =  \{ (0,1), (0,-1) \}$  and $\{ f_2 = 0 \} \cap \{ x_1 = 0 \} =  \{ (1,0), (-1,0) \}$.
Moreover  $\{ f_1 = 0 \} \cap \{ f_2 = 0 \}  =  \{ (\pm\frac{1}{\sqrt{3}},\pm\frac{1}{\sqrt{3}}) \} $.  
\bigskip

The theorem \ref{asystab} may be extended to prove the global asymptotic stability of an equilibrium point.

\begin{corollary} \label{stabglob}  If $\Omega = \R^n$, $(a_i,b_i) = R$ for $i=1,\dots,n$,  and the hypotheses of theorem \ref{asystab} hold, then $O$ is globally asymptotically stable for (\ref{sysn}).
\end{corollary}
{\it Proof.} Under the above hypotheses, the Liapunov function $V$ of the theorem \ref{asystab} is definied on all of $\R^n$, with $\dot V \leq 0$ on all of $\R^n$.  By LaSalle invariance principle, every solution $\gamma$ tends to the largest invariant set $E$ contained in $\dot V = 0$.  As proved in theorem \ref{asystab}, $E = \{O\}$. This gives the global attractivity of the origin, hence its global  asymptotic stability.
\hfill$\clubsuit$

As a consequence of theorem \ref{asystab}, we prove an existence result for periodic solutions of small periodic perturbations of the equation  \ref{equaliev}.

\begin{corollary} \label{persol}  
If the hypotheses of theorem \ref{asystab} hold, and $h(t,X,X',\varepsilon)  \in C^1(\R \times \Omega \times \R^{n} \times [0,\overline{\varepsilon}),\R^n)$, $\overline{\varepsilon} > 0$, then there exists $\varepsilon^* \in(0,\overline{\varepsilon})$ such that the equation  (\ref{equalievper}) has a periodic solution $\phi(t, \varepsilon )$ for all $\varepsilon \in (0, \varepsilon^*)$, such that $\phi(t, \varepsilon )$ tends to the null solution as $\varepsilon \rightarrow 0$.
\end{corollary}
{\it Proof.}  
The asymptotic stability of a critical point of an autonomous system is equivalent to its uniform asymptotic stability, hence we may apply theorem 15.9 in \cite{Y} in order to prove the statement.
\hfill$\clubsuit$

\end{document}